\documentclass[12pt]{amsart}
\usepackage{graphics,color,epsfig}
\usepackage{amsfonts}
\usepackage{amssymb}

\newcommand{\C}{\mathbb{C}}

\newcommand{\Q}{\mathbb{Q}}
\newcommand{\R}{\mathbb{R}}
\newcommand{\Z}{\mathbb{Z}}

\newcommand{\CCC}{\mathcal{C}}
\newcommand{\HHH}{\mathcal{H}}
\newcommand{\MMM}{\mathcal{M}}

\newcommand{\cs}{\sigma}

\DeclareMathOperator{\area}{area}
\DeclareMathOperator{\hol}{hol}
\DeclareMathOperator{\SL}{SL}
\DeclareMathOperator{\GL}{GL}

\newtheorem{thm}{Theorem}[section]
\newtheorem{lemma}[thm]{Lemma}
\newtheorem{prop}[thm]{Proposition}
\newtheorem{coro}[thm]{Corollary}

\newtheorem{definition}[thm]{Definition}

\begin{document}

\title[Topological Dichotomy and Strict Ergodicity]
{Topological dichotomy and strict ergodicity for translation surfaces}

\author{Yitwah Cheung, Pascal Hubert, Howard Masur}

\address{San Francisco State University \\ 
San Francisco, CA 94132, U.S.A.}
\email{cheung@math.sfsu.edu}

\address{LATP, Case cour A, Facult\'e des Sciences de Saint J\'er\^ome, \\
Avenue Escadrille Normandie-Niemen \\
13397, Marseille Cedex 20, France}
\email{hubert@cmi.univ-mrs.fr}

\address{University of Illinois at Chicago \\
Chicago, Illinois 60607, U.S.A.}
\email{masur@math.uic.edu}

\thanks{The first author acknowledges the University of Illinois at Chicago for financial support when this work was done.}  
\thanks{The second author acknowledges the University of Chicago and the University of Illinois at Chicago for hospitality and financial support when this work was done.}    
\thanks{The research of the third author is partially supported by NSF grant DMS0244472.}

\subjclass{32G15, 30F30, 30F60, 58F18} \keywords{Abelian differentials, Teichm\"uller disc, unique ergodicity, minimality}
\date{\today}

\begin{abstract}
In this paper the authors find examples of translation surfaces 
  that have infinitely generated Veech groups, satisfy the 
  topological dichotomy property that for every direction either 
  the flow in that direction is completely periodic or minimal,
  and yet have minimal but non uniquely ergodic directions. 
\end{abstract}
\maketitle

\section{Introduction and Statement of Theorems}
Let $(X,\omega)$ be a translation surface where $X$ is a closed surface of genus at least $2$.  Equivalently, $X$ is a Riemann surface, and $\omega$ is a holomorphic $1$-form on $X$.  For each $\theta\in[0,2\pi)$ there is a vector field defined on the complement of the zeroes of $\omega$ such that $\arg\omega=\theta$ along this vector field.  The corresponding flow lines are denoted $\phi_\theta$.  A great deal of work has been done to try to understand the dynamics of $\phi_\theta$. For a countable set of $\theta$ there is a flow line of $\phi_\theta$ joining a pair of zeroes of $\omega$.  These flow lines are called saddle connections.  For any $\theta$ such that there is no saddle connection in direction $\theta$, it is well-known that the flow is minimal.  
  Veech (\cite{Ve2}) introduced an important class of translation surfaces, now called \emph{Veech surfaces}.  
They are defined by the property that the stabilizer $\SL(X,\omega)$ of the surface under the action of $\SL(2,\R)$ on the moduli space of translation surfaces (or abelian differentials) is a lattice.\footnote{The term \emph{lattice surface} is also commonly used to describe these surfaces.}  
These surfaces satisfy (\cite{Ve2}) optimal dynamics in that for every direction $\theta$, either the flow lines of $\phi_\theta$ are all closed or saddle connections and the surface decomposes into a union of cylinders, or the flow is minimal and uniquely ergodic.   

On the other hand Veech \cite{Ve1} had previously found examples of skew rotations over the circle that are not strictly ergodic.  Namely, the orbits are dense but not uniformly distributed.  Veech's examples can be interpreted (\cite{MaTa}) in terms of flows on $(X,\omega)$ where $X$ has genus $2$ and $\omega$ has a pair of simple zeroes.  
%Take two copies of the standard torus $\R^2/\Z^2$ and mark off a segment along the vertical axis from $(0,0)$ to $(0,\alpha)$, where $0<\alpha<1$.  Cut each torus along the segment and glue pairwise along the slits.  The resulting surface $(X_\alpha,\omega)$ is the connected sum of the pair of tori and is a branched double cover over the standard torus, branched over $(0,0)$ and $(0,\alpha)$.  These two endpoints of the slits become the zeroes of order one of $\omega$.  There are a pair of circles on $X_\alpha$ such that the first return map of $\phi_\theta$ to these circles gives a skew rotation over the circle.  If $\alpha$ is irrational, then there are directions $\theta$ such that the flow $\phi_\theta$ is minimal but not uniquely ergodic.  These reproduce the original Veech examples. 
Veech's theorem together with his examples raised the issue of what can be said about the dynamics of flows on surfaces $(X,\omega)$ that are not Veech surfaces, both from the measure-theoretic and topological points of view.  

We say that $(X,\omega)$ satisfies {\em topological dichotomy}, if for every direction, either the flow is minimal, or every flow line is closed or a saddle connection.  This is equivalent to saying that if there is a saddle connection in direction $\theta$, then there is a cylinder decomposition of the surface in that direction. 
A translation surface $(X,\omega)$ is said to satisfy {\em strict ergodicity} if every minimal direction is uniquely ergodic.  Veech surfaces satisfy both topological dichotomy and strict ergodicity.  

In genus $2$ McMullen (\cite{Mc4}) showed that {\em every} surface that is not a Veech surface does not satisfy topological dichotomy.  Cheung and Masur (\cite{ChMa}) showed that {\em every} surface of genus $2$ that is not a Veech surface does not satisfy strict ergodicity.  Thus, neither generalisation of the notion of a Veech surface can be realized in genus $2$.  

From the measure theoretic point of view, in the paper \cite{MaSm} it was shown that a generic translation surface in any moduli space does not satisfy the strict ergodicity property.   

In this paper we prove the following theorem.
\begin{thm}
\label{thm:main}
There are examples of translation surfaces that satisfy topological dichotomy but not strict ergodicity.  Moreover, these examples have infinitely generated $\SL(X,\omega)$.  
\end{thm}
These examples are based on a construction of Hubert and Schmidt (\cite{HuSc}) whose original motivation was to give examples of translation surfaces with infinitely generated $\SL(X,\omega)$.  They are double covers of genus $2$ Veech surfaces that are branched over a singularity and a aperiodic connection point.  Recall that a \emph{connection point} is defined by the property that every segment from a singularity to the point extends to a saddle connection.  It is \emph{periodic} if its orbit under $\SL(X,\omega)$ is finite.  The existence of aperiodic connection points was first established in \cite{HuSc} for any non-arithmetic Veech surface in the stratum $\HHH(2)$.  

We also give a billiard example:

\begin{thm}\label{thm:billiard}
Let $\Delta$ be the rational billiard table of the 
  $(\frac{2\pi}{5},\frac{3\pi}{10},\frac{3\pi}{10})$ triangle and 
$(\Hat X,\Hat\omega)$ the associated translation surface obtained by the unfolding process.
Then $(\Hat X,\Hat\omega)$ satisfies topological dichotomy but not strict ergodicity.  
\end{thm}

A slightly weaker condition than topological dichotomy is {\em complete periodicity.}  A direction $\theta$ is said to be completely periodic if every flow line in direction $\theta$ is closed.  In that case, the surface decomposes into cylinders.  Each cylinder is maximal in the sense that it cannot be enlarged.  The boundary of these cylinders consists of saddle connections.  The surface $(X,\omega)$ is said to be {\em completely periodic} if any direction that has {\em some} cylinder in that direction is completely periodic.  Calta (\cite{Ca}) gave examples of families of genus $2$ translation surfaces that satisfy complete periodicity and yet are not lattice surfaces.  These examples were independently discovered by McMullen \cite{Mc4} who further showed that they do not satisfy topological dichotomy.  

It is worth noting that if a surface is a finite branched cover over a Veech surface, then it trivially satisfies complete periodicity.  Namely, any direction with a cylinder on the cover gives rise to a completely periodic direction on the Veech surface, which in turn means that the direction is completely periodic on the cover.  If the cover is branched over a set of connection points, at most one of which is aperiodic, then it satisfies topological dichotomy as well.  

Smillie and Weiss (\cite{SmWe}) have shown that a surface satisfying both topological dichotomy and strict ergodicity need not be a Veech surface.  These examples are also based on the construction of Hubert and Schmidt (\cite{HuSc}) but differ from our examples in that they are branched over a single aperiodic connection point.  Smillie and Weiss's construction proves strict ergodicity for every cover over a Veech surface branched over one point. If the branched point is not a connection point, then the cover does not satisfy the topological dichotomy. This means that there are also examples that satisfy strict ergodicity and not topological dichotomy.

We actually prove the following more general theorem.  
\begin{thm}
\label{thm:general}
Let $(\hat X,\hat \omega)$ be a branched double cover over any lattice surface $(X,\omega)$ branched over the singularity and a regular aperiodic point.   Then $(\hat X, \hat\omega)$ does not satisfy strict ergodicity.
\end{thm}

The classification of Veech surfaces in higher genus is far from complete.  Theorem \ref{thm:general} does apply to infinitely many surfaces $X$ of genus higher than 2, such as pairs of regular $n$-gons with parallel sides identified, the original examples of Veech surfaces (see \cite{Ve2}). McMullen proved the existence of infinite families of Veech surfaces in genus 3 and 4 (see \cite{Mc6}). Bouw and M\"oller gave another family of examples obtained by  deep algebraic methods (see \cite{BoMo}). Nevertheless Theorem \ref{thm:main} only gives examples of some coverings of {\it genus 2} Veech surfaces.  The existence of aperiodic connection points is an open question when the genus of the Veech surface is greater than 2.  Thus, topological dichotomy for coverings of Veech surfaces of genus greater than 2 is not known.  

%The regular $2n$-gon with opposite sides identified defines a Veech surface in the hyper-elliptic component of $\HHH(2g-2)$ if $n$ is even, where $n=2g$ and in $\HHH(g-1,g-1)$ if $n$ is odd, where $n=2g+1$.  Our techniques can also be applied to give results for some branched covers of these surfaces.  
%\begin{thm}
%\label{thm:2n-gons}
%Let $D$ be the $\SL(2,\R)$ orbit of the regular $2n$-gon and set $\beta=(4g-3,1)$ if $n$ is even and $\beta=(2g-1,g-1,g-1,1)$ if $n$ is odd.  Let $\MMM_D(\beta)$ be the moduli space of translation surfaces in $\HHH(\beta)$ that are branched double covers of some surface in $D$.  (The choice of $\beta$ guarantees $\MMM_D(\beta)$ is nonempty.)  Then any surface in $\MMM_D(\beta)$ that is not a Veech surface does not satisfy strict ergodicity either.  
%\end{thm}

\subsection*{Acknowledgments}
The authors would like to thank Matt Bainbridge, Emanuel Nipper and John Smillie for useful comments.

\section{Background}
In this section, we recall some results that we will use later in the proof of the main theorem.  

\subsection{Veech surfaces}
Let $(X,\omega)$ be a translation surface.  The stabilizer of $(X,\omega)$ under the $\SL(2,\R)$ action is called the {\it Veech group} of $(X,\omega)$ and is denoted by $\SL(X,\omega)$.  
A more intrinsic definition is the following.  An {\it affine diffeomorphism} is an orientation preserving homeomorphism of $X$ which is affine in the charts of $\omega$ and permutes the zeroes of $\omega$.  The derivative (in the charts of $\omega$) of an affine diffeomorphism defines an element of $\SL(X,\omega)$ and vice versa.  

A translation surface $(X,\omega)$ is a {\it Veech surface} if its Veech group is a lattice in $\SL(2,\R)$.

We recall that a Veech surface satisfies both topological dichotomy and strict ergodicity by a theorem of Veech \cite{Ve2}.  (This is sometimes called Veech dichotomy in the literature.)  Moreover a periodic direction on a Veech surface $(X,\omega)$ is fixed by a parabolic affine diffeomorphism.  As $\SL(X,\omega)\backslash \SL(2,\R)$ has a finite number of cusps there are a finite number of periodic directions up to the action of the affine group of diffeomorphisms.  Therefore, there are a finite number of cylinder decompositions; each one is associated to a cusp. 

A Veech surface is {\it arithmetic} if it is a covering of a flat torus, ramified over one point;  otherwise, it is a {\it non-arithmetic} Veech surface.

\subsection{Classification of Veech surfaces in genus 2}
There are two moduli spaces or strata of abelian differentials in genus 2: $\HHH(1,1)$ the stratum of differentials with two zeroes of order one and $\HHH(2)$ the stratum of differentials with a single zero of order $2$.  

McMullen gave a complete classification of non-arithmetic Veech surfaces in genus 2 \cite{Mc1, Mc2, Mc3, Mc4, Mc5} (see also Calta \cite{Ca} for some partial results and a different approach).
McMullen proved that in $\HHH(2)$ there are infinitely many non-arithmetic Veech surfaces and classified them.  Each of them is defined over a real quadratic field which means that, up to normalization, all the relative periods lie in a quadratic field. 

Every $\SL(2,\R)$ orbit of a Veech surface in $\HHH(2)$ contains a $L$-shaped surface and we have 
\begin{thm} [McMullen \cite{Mc3}]
  The $L$-shaped surface $L(a,b)$ is a Veech surface if and only if 
 $a$ and $b$ are rational or $a = x + z \sqrt{d}$ and 
 $b = y + z \sqrt{d}$ for some $x, y, z \in \Q$ with 
 $x + y = 1$ and $d \geq 0$ in $\Z$.
\end{thm}
\noindent (see figure  \ref{L(a,b)} for the definition of $L(a,b)$).

\begin{figure}[ht] 
\includegraphics{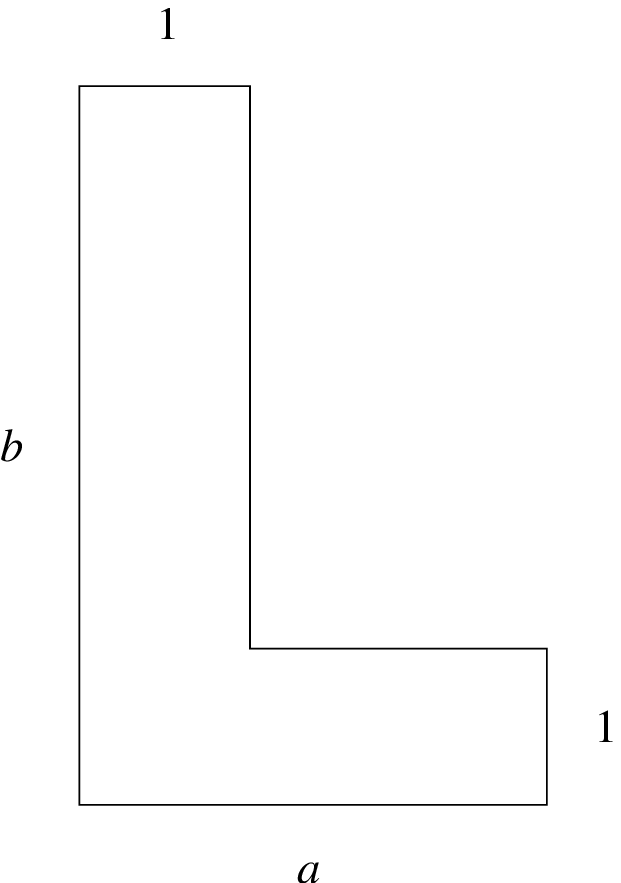}  
\caption{$L(a,b)$}
\label{L(a,b)}
\end{figure}

The situation is very different in $\HHH(1,1)$.  There is only one non-arithmetic Veech surface (up to action of $\SL(2,\R)$): the surface obtained from the regular decagon by gluing opposite sides together (see \cite{Mc4, Moe, Mc5}.

\subsection{Periodic directions in $\HHH(2)$} \label{periodic-directions-inH(2)}
In a completely periodic direction, a surface $(X,\omega)$ in $\HHH(2)$ decomposes into one or two cylinders (\cite{Zo}).  If $(X,\omega)$ is a non-arithmetic Veech surface, every periodic direction induces a two cylinders decomposition \cite{Ca}.  
For simplicity, assume that the horizontal direction is periodic.  There are three horizontal saddle connections.  Two of them are exchanged by the hyperelliptic involution.  We will call them {\it non Weierstrass saddle connections}.  The third one contains a Weierstrass point and is fixed by the hyperelliptic involution; it will be called a {\it Weierstrass saddle connection}.  

\begin{figure}[h] 
 \includegraphics{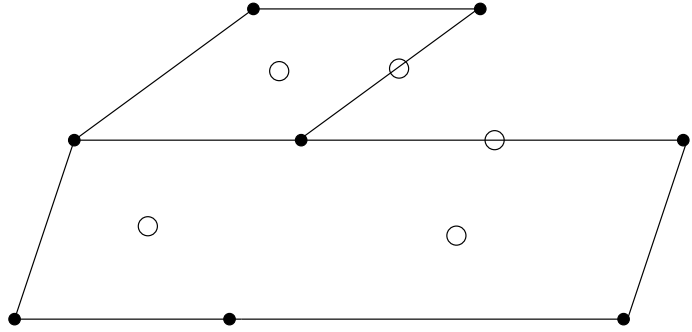}  
 \caption{Periodic direction in $\HHH(2)$. Black points correspond to the singular point and white points to the five other Weierstrass points. Parallel sides of the same length are identified.}
 \label{fig:standard-picture}
\end{figure}

The boundary of one cylinder consists of the two saddle connections exchanged by the hyperelliptic involution.  We call it the {\it non self-gluing cylinder}.  The other one has two saddle connections on each component of its boundary.  One of them contains a Weierstrass point and lies on both components.  It is the {\it self-gluing} cylinder.  We summarize this in figure \ref{fig:standard-picture}.

\subsection{Infinitely generated Veech groups and connection points}
Constructions of translation surfaces with infinitely generated Veech groups can be found in \cite{HuSc} and \cite{Mc2}.  We recall some ingredients in Hubert-Schmidt's approach.  

%{\it Definitions:} A point $P$ on a Veech surface is {\it periodic} if its orbit under the affine group is finite.  
%A point $P$ on a translation surface is a {\it connection point} if every separatrix containing $P$ is a saddle connection.  

\begin{thm} [Hubert-Schmidt \cite{HuSc}]
Any double covering of a Veech surface ramified over a singularity and a aperiodic connection point has infinitely generated Veech group and satisfies topological dichotomy.  
\end{thm}

M\"oller (\cite{Moe}) proved that the Weierstrass points are the only periodic points on a non-arithmetic Veech surface of genus 2.  Following Calta and McMullen, one can show that a point on the $L$-shaped Veech surface $L(a,b)$, with $a$, $b \in \Q(\sqrt{d})$, is a connection point if and only if its coordinates belong to $\Q(\sqrt{d})$ (the origin is the singular point).  Consequently, if a point in $L(a,b)$ has a coordinate which is not in $\Q(\sqrt{d})$ it is contained in at most one saddle connection (see \cite{HuSc}).

\subsection{Space of coverings}
Given a translation surface $(X,\omega)$ with singular set $\Sigma$ one can consider the moduli space $\MMM$ of all double covers $(\hat X,\hat \omega)$ branched over two points, one in $\Sigma$ and the other not in $\Sigma$, up to biholomorphic equivalence.  An explicit way to realize an element of $\MMM$ is to take two copies of $(X,\omega)$ each marked with a segment $\cs$ from some point $P_0\in\Sigma$ to a nearby point $P_1\not\in\Sigma$.  Cut both copies along $\cs$ and glue them pairwise.  We call $\cs$ a {\em slit} and refer to this construction as a {\em slit construction}.  Note that Smillie and Weiss's examples are not formed in this manner.  

\begin{lemma}
\label{lem:connected}
%The number of connected components of $\MMM$ is bounded by the number of orbits for the action of $\SL(X,\omega)$ on $\Sigma$.  
Each surface in $\MMM$ can be obtained via a slit construction.  Moreover, two surfaces lie in the same connected component if they are branched over zeroes that correspond under the $\SL(2,\R)$ action.  
\end{lemma}
\begin{proof}
Let $\pi:(\hat X,\hat \omega)\to(X,\omega)$ be the natural projection and $\tau:\hat X\to \hat X$ the involution such that $\pi\circ\tau=\pi$.  
Let $\hat P_0,\hat P_1\in\hat X$ be the branching points lying over $P_0$ and $P_1$, respectively.  
Let $\hat\gamma$ be any separating multi-curve on $\hat X$ such that $\tau(\hat\gamma)$ is homotopic to $\hat\gamma$.  (By a multi-curve we mean a union of disjoint simple closed curves.  To construct $\hat\gamma$, choose any simply-connected subset of $X$ whose complement is a finite union of arcs containing $P_0$ and $P_1$; the closure of either one of its two lifts is a sub-surface whose boundary is fixed by $\tau$.)  It is a standard fact that there is a geodesic with respect to the flat structure of $\hat\omega$ in the class of $\hat\gamma$, again denoted by $\hat\gamma$.  It is uniquely determined if and only if no component of $\hat\gamma$ is homotopic to the core curve of a cylinder in $\hat X$ and it has no transverse self-intersections.  We may choose it so that $\tau(\hat\gamma)=\hat\gamma$.  Let $\hat X_1,\hat X_2$ be the two components in the complement of $\hat\gamma$ that are interchanged by $\tau$.  

It is clear that both $\hat P_0,\hat P_1\in\hat\gamma$; otherwise, one of them would lie in the interior of either $\hat X_1$ or $\hat X_2$, and since $\tau$ fixes the point, it would not interchange $\hat X_1$ and $\hat X_2$.  Let $\cs_1$ be a saddle connection in $\hat\gamma$ with an endpoint at $\hat P_1$.  Note that $\tau(\cs_1)\neq\cs_1$ because the interior of $\cs_1$ would contain a point fixed by $\tau$, which is impossible.  The pair $\cs_1\cup\tau(\cs_1)$ divides a neighborhood of $\hat P_1$ into two components interchanged by $\tau$ so that there are an equal number $n$ of segments of $\hat\gamma$ on each side of $\cs_1\cup\tau(\cs_1)$.  Since $\hat\gamma$ is separating, $n$ is even.  The total angle at $\hat P_1$ is $4\pi$ so that there is a total angle of $2\pi$ on each side of $\cs_1\cup\tau(\cs_1)$.  Since $\hat\gamma$ is a geodesic, the angle between incoming and outgoing segments is at least $\pi$.  This is not possible for if $n\geq 2$ so we must have $n=0$.  Hence, there is a unique component of $\hat\gamma$ that contains $\hat P_1$.  This component is fixed by $\tau$ and also contains $\hat P_0$.  

%We now claim that there are in fact no other components of $\hat\gamma$ fixed by $\tau$.  Indeed, otherwise such a component must consist of in order, a concatenation of saddle connections $$\cs_0*\ldots *\cs_0'*\tau(\cs_0')*\ldots*\tau(\cs_0),$$ where $\cs_0$ and $\tau(\cs_0')$ begin at $\hat P_0$, while $\cs_0'$ and $\tau(\cs_0)$ end at $\hat P_0$ (possibly $\cs_0=\cs_0'$).  Since $\cs_0$ and $\tau(\cs_0)$ subdivide the total angle at $\hat P_0$ into equal pieces, and the same is true of $\cs_0'$ and $\tau(\cs_0')$, this component of $\hat\gamma$ has self-intersections, which is a contradiction, proving the claim.  Hence, there is exactly one component of $\hat\gamma$ fixed by $\tau$, namely the one containing $\hat P_1$, while all remaining components, if any, are interchanged by $\tau$ in pairs.  

The projection $\pi(\hat\gamma)$ of $\hat\gamma$ to $X$ consists of an arc joining $P_1$ to $P_0$ together with a finite union of loops.  A closed curve in $X\setminus \{P_0,P_1\}$ lifts to a closed curve if and only if it intersects $\pi(\hat\gamma)$ an even number of times.  Let $\gamma$ be a curve in $X$ joining $P_1$ to $P_0$ that is homologous to $\pi(\hat\gamma)$ as chains relative to $\Sigma\cup P_1$.  It is easy to see that we may represent $\gamma$ by an closed embedded arc.  Since $\gamma$ and $\pi(\hat\gamma)$ are homologous (mod $2$), the surface $(\hat X,\hat\omega)$ can be obtained via a slit construction using  $\gamma$.  Notice here we are not assuming that $\gamma$ is a geodesic arc. 

Now fixing the one endpoint $P_0$, and moving the other endpoint, we may continuously deform $\gamma$ to a fixed geodesic segment from $P_0$ to a point in a small neighborhood of $P_0$.  The family of arcs along this path all determine double covers by a slit construction and thus all belong to the same connected component of $\MMM$.  Thus any two double covers branched over $P_0$ all belong to the same component.  %This shows that  number of connected components is bounded by the order of $\Sigma$.  If an element of $\SL(X,\omega)$ takes $P_0$ to $P_0'$ then the corresponding connected components are the same, because there is a double cover of $(X,\omega)$ that can be realized via the slit construction in two ways, using a slit at either $P_0$ or at $P_0'$.  
\end{proof}

Let $\HHH(\beta)$ be a stratum of the moduli space of abelian differentials of genus $g$ and $D$ be the $\SL(2,\R)$-orbit of a Veech surface $(X,\omega)$.  The moduli space $\MMM_D(\beta)\subset\HHH(\beta)$ of all possible coverings of surfaces in $D$ is a closed $\SL(2,\R)$-invariant orbifold.  Eskin-Marklof-Morris \cite{EsMaMo} proved an analogue of Ratner's theorem for the action of the horocycle flow on $\MMM_D(\beta)$.  They classified the invariant measures and the orbit closures.  We state a very weak version of their result which is enough for our purpose.  
\begin{thm} [Eskin-Marklof-Morris \cite{EsMaMo}]
\label{thm:dense}
Assume $(X,\omega)$ is a Veech surface and suppose that $(\Hat X,\Hat\omega)\in\HHH(\beta)$ is obtained from $(X,\omega)$ by a slit construction.  Then the $\SL(2,\R)$-orbit of any translation surface in $\MMM_D(\beta)$ that is not a Veech surface is dense in the connected component of $\MMM_D(\beta)$ that contains it.  
\end{thm}

%Let $\XXX$ be the fiber bundle over $\SL(2,\R)X$ whose fiber over $gX$ is the surface $g(X\setminus\text{zeroes}(\omega))$ for all $g\in\SL(2,\R)$.  A point in $\XXX$ is a pair $(X',q)$ where $X'$ belongs to $\SL(2,\R)X$ and $q$ is a regular point in $X'$.  It can be shown that $\MMM_D(\beta)$ is a finite branched cover of $\XXX$.  

\subsection{Non uniquely ergodic directions}
Let $(X,\omega)$ be a translation surface, $\gamma$ be a dividing geodesic curve on $X$ such that $X \setminus \gamma$ has two connected components $\Omega_1$, $\Omega_2$.  The triple $(\Omega_1, \Omega_2,\gamma)$ is a {\it splitting} of $(X,\omega)$. 

The following theorem gives a criterion to get nonuniquely ergodic directions:

\begin{thm} [Masur-Smillie \cite{MaSm}]
\label{thm:criterion}
Let $(\Omega_1^n,\Omega_2^n,\gamma_n)$ be a sequence of splittings of $(X,\omega)$ and assume the directions of the holonomy vectors $\gamma_n$ converge to some direction $\theta$.  Let $h_n$ be the component of $\gamma_n$ in the direction perpendicular to $\theta$ and $a_n = \area(\Omega_1^n \Delta \Omega_1^{n+1})$, the area of the regions exchanged between consecutive splittings. If
 \begin{enumerate}
 \item $\displaystyle\sum_{n= 1}^\infty a_n < \infty,$
 \item there exists $c >0$ such that $\area(\Omega_1^n) > c$, $\area(\Omega_2^n) > c$ for all $n$, and
 \item $\displaystyle\lim_{n \to \infty} h_n = 0,$ 
 \end{enumerate}
 then $\theta$ is a nonergodic direction.
\end{thm}

In the sequel, we will construct uncountably many sequences of splittings satisfying the above three conditions with distinct limiting directions $\theta$.  As there are only countably many nonminimal directions, there must be uncountably many minimal and non uniquely ergodic directions.

\section{Preliminary lemmas}

\subsection{Parallelograms and change of area}
Let $(X,\omega)$ be a Veech surface, $\cs$ a slit on $X$ joining a singularity $P_0$ to another point $P$.  Let $(\hat X,\hat\omega)$ be the cover of $X$ obtained by gluing two copies of $X$ together along $\cs$.  We will call $(\hat X,\hat\omega)$ the {\it cover induced by the slit} $\cs$.  Let $\hat \cs_1,\hat \cs_2$ the two copies of $\cs$ on $\hat X$.  Let $\hat \cs=\hat \cs_1-\hat \cs_2$ the separating curve on $\hat X$ which divides $\hat X$ into components $\Omega_1,\Omega_2$.  Suppose $\cs$ is contained in a cylinder $\CCC'$ such that $P_0$ is on one component of the boundary of $\CCC'$ and $P$ is on the other.  Note we are not assuming $\CCC'$ is a maximal cylinder--the boundary component that contains $P$ will not contain any singularities.  Suppose $\cs'$ is another segment of constant slope in $\CCC'$ also joining $P_0$ to $P$.  Assume $\cs$ and $\cs'$ intersect an odd number of times in their interior. 
\begin{lemma}
\label{lemma:area}
The two lifts $\hat \cs_1',\hat \cs_2'$ of $\cs'$ also divide $\hat X$ into two components $\Omega'_1,\Omega_2'$ and either $\area(\Omega_1\Delta \Omega_1')\leq 2\area(\CCC')$ or $\area(\Omega_1\Delta\Omega_2')\leq 2\area(\CCC')$.
\end{lemma}
\begin{proof}
Let $m$ be the number of intersections.  We have  $\cs-\cs'$ is homologous to $(m+1)\beta$, where $\beta$ is the core curve of $\CCC'$.  Since $m$ is odd, d $m+1$ is even, and so the lift of $(m+1)\beta$ is $\frac{m+1}{2}\hat\beta$, where $\hat\beta$ is a closed curve.  Then $\hat \cs_1$ is homologous to $\hat\cs_1'+\frac{m+1}{2}\hat \beta$ and $\hat \cs_2$ is homologous to $\hat \cs_2'+\frac{m+1}{2}\hat\beta$.  Since $\hat \cs_1$ is homologous to $\hat \cs_2$, we have $\hat \cs_1'$ homologous to $\hat \cs_2'$ which implies the first statement. 

We prove the second statement.  The intersections of $\cs'$ with $\cs$ divides each into $m+1$ segments of equal length.  There are $\frac{m+1}{2}$ congruent parallelograms $P_j$ in $\CCC'\setminus(\cs\cup \cs')$ such that 
\begin{itemize}
  \item each segment of $\cs$ and each segment of $\cs'$ lies on the boundary of exactly one parallelogram (see figure \ref{figure:area}).  
  \item a point in any  parallelogram can be joined to a point in any other by a path that crosses $\cs\cup \cs'$ an even number of times. 
\end{itemize}

\begin{figure}[ht] 
 \includegraphics{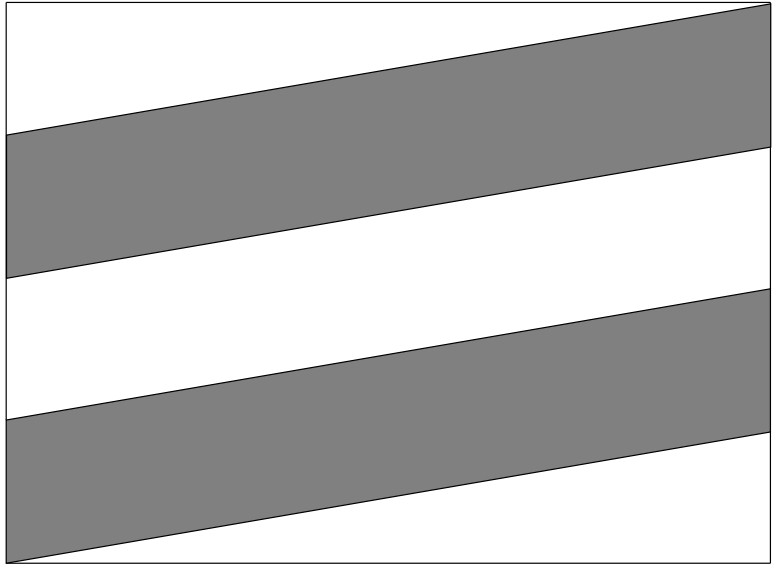}  
 \caption{change of area}
 \label{figure:area}
\end{figure}

Lift these parallelograms to parallelograms $\hat P_i$ in $\hat X$.  The conditions mean that all parallelograms belong to $(\Omega_1\cap \Omega_1')\cup (\Omega_2\cap \Omega_2')$ or they all belong to $(\Omega_1\cap \Omega_2')\cup(\Omega_2\cap\Omega_1')$.  That is to say, in the first case they belong to $\Omega_1\Delta \Omega_2'$, and in the second to $\Omega_1\Delta\Omega_1'$.  Every other point of $\hat X$ belongs to $\Omega_1\Delta\Omega_2$ in the first case and $\Omega_1\Delta\Omega_2'$ in the second.  Since the sum of the areas of the parallelograms on $X$ is at most the area of $\CCC'$, we have the result.  
\end{proof}

\begin{definition}
We will denote by $a(\cs,\cs')$, the minimum of the two areas in the conclusion of the last lemma and call it the change in area.  
\end{definition}

\begin{definition}
If $\cs$ is a slit contained in a saddle connection $\gamma$ then the {\it ratio} of $\cs$ is defined to be the length of $\cs$ divided by the length of $\gamma$.  If $\cs$ is not contained in a saddle connection its ratio is taken to be zero.  The ratio of $\cs$ will be denoted by $\rho(\cs)$.  
\end{definition}

\begin{lemma} \label{ratio-cross-product}
Let $(X,\omega)$ be a Veech surface of area one, $\cs$ a slit on $(X,\omega)$ that lies on a saddle connection $\gamma$, and $(\hat X,\hat\omega)$ the induced covering by $\cs$.  
We assume that the horizontal direction is a periodic direction and that $\gamma$ crosses $k$ horizontal cylinders each of them at most $\ell$ times.  We also assume that $\cs$ is contained in a maximal horizontal cylinder $\CCC_1$.  Let $\CCC_1'$ be the smallest horizontal subcylinder of $\CCC$ containing $\cs$ and $T$ the affine Dehn twist in the cylinder $\CCC_1'$.  Let $\cs_n$ be the image of $\cs$ under $T^{2n}$.  Then the change of area $a(\cs,\cs_n)$ satisfies $a(\cs,\cs_n) \leq C(\ell)\rho(\cs)$ where $C(\ell)$ is a constant depending only on $\ell$.
\end{lemma}
\begin{proof}
The change of area is bounded by twice the area of $\CCC_1'$, by Lemma~\ref{lemma:area}.  Let $\CCC_1, \dots, \CCC_k$ be the horizontal cylinders that $\gamma$ crosses, $h_1, \dots, h_k$ are the heights of these cylinders.  Let $$\kappa = \frac{\max\{h_i, \ 1 \leq i \leq k\}}{\min\{h_j, \ 1 \leq i \leq k\}}.$$
 Let $y$ be the vertical component of $\hol(\cs)$.  We have $$\area(\CCC_1') = \frac{y}{h_1}\times \area(\CCC_1)\leq \frac{y}{h_1}.$$  For $i = 1 \mbox{ to } k$, let $a_i$ be the number of times that $\gamma$ crosses $\CCC_i$.  We recall that $a_i \leq \ell$. The vertical component of $\hol(\gamma)$ is $a_1 h_1 +\cdots +a_k h_k$.  Thus 
$$\area(\CCC_1') \leq \frac{y}{a_1 h_1 +\cdots +a_k h_k}\times 
\frac{a_1 h_1 +\cdots +a_k h_k}{h_1} \leq \rho \  k \ell \kappa.$$
The number $\kappa$ is bounded because, as $(X,\omega)$ is a Veech surface, there are finitely many cusps and, of course, finitely many cylinders in each cylinder decomposition.
Thus $\kappa$ is bounded by a constant (depending only on the geometry of the surface $X$).  Therefore there exists $C(\ell)$ such that $a(\cs,\cs_n) \leq  C(\ell) \rho(\cs)$.
\end{proof}

\subsection{Irrationality and small ratios}
Lemma~\ref{ratio-cross-product} says that the change in area is bounded in terms of the ratio of the slit in the case that the slit is part of a saddle connection.  Thus it will be important to find sequences of slits with small ratio. Lemma~\ref{keylem:twist} below  gives  conditions under which this is possible.  Proposition~\ref{prop:angle} will allow us to find a cylinder to which Lemma~\ref{keylem:twist} can be applied.  We will use these ideas in the proof of Theorem~\ref{thm:general} where we will show that we can make the new ratios small in all situations. 

\begin{lemma}
\label{keylem:twist}
Let $\cs$ be a slit contained in maximal horizontal cylinder $C$ and suppose that the vertical component $v_\cs$ of $\cs$ is an irrational multiple of the height of $C$.  Then, given any $\delta>0$, there exists an even $n$ such that the $n$th twist $\cs_n$ has ratio less than $\delta$.  
\end{lemma}
\begin{proof}
Consider the set $S$ of holonomy vectors of all saddle connections or loops whose vertical components are bounded by $\frac{v_\cs}{\delta}$.  This set $S$ is invariant under the parabolic subgroup $U$ that fixes the horizontal direction.  We claim that $S$ is a finite union of orbits of $U$.  

The generator of $U$ is given by the matrix 
$$
\left(\begin{array}{cc}
             1  &   \alpha \\
             0  &   1
\end{array}\right)
$$ 
where, in the horizontal decomposition into cylinders $C_i$, there are integers $n_i$ such that $\alpha=n_iC_i$ independent of $i$.  The generator acts on inverse slopes by translating the fixed amount $\alpha$.  That means every orbit contains a representive whose inverse slope is between $0$ and $\alpha$.  Since the vertical coordinate of the representative is bounded above, so is the horizontal coordinate.  That means that the holonomy vector of the representative is bounded.  There are only finitely many saddle connections with bounded holonomy, which implies there are only finitely many orbits under this action of $U$.  

On the other hand we twist $\cs$ in a subcylinder whose height is an irrational multiple of the height $C$ and whose circumference is the same as that of $C$.  The inverse slope of the resulting $\cs_n$ is $$\frac{x}{y}+n\lambda\alpha$$ where $\frac{x}{y}$ is the inverse slope of $\cs$ and $\lambda$ is irrational.  A pair $\cs_n,\cs_{n'}$ cannot be in the same $U$-orbit since the equation $(n-n')\lambda=m$ does not have a solution for $n\neq n'$.  This implies there is some twist $\cs_n$ with $n$ even so that if $\cs_n$ is part of a saddle connection, then the vertical component $v'$ of this saddle connection is greater than $\frac{v_\cs}{\delta}$ and therefore the ratio $\frac{v_\cs}{v'}<\delta$.  
\end{proof}

\begin{definition}
We say two cylinders $C_1$ and $C_2$ in different directions are \emph{rationally related} if, after normalizing by the action of $\SL(2,\R)$ so that one is horizontal and the other is vertical, the ratio of the height of $C_1$ to the circumference of $C_2$ and the ratio of the height of $C_2$ to the circumference of $C_1$ are both rational.  
\end{definition}

\begin{definition}
Let $C$ be a cylinder and $P$ a singularity on the boundary of $C$.  We define $\alpha(P,C)$ to be the open set of unit tangent vectors $v$ at $P$ such that an initial segment of $v$ lies in the interior of $C$.  
\end{definition}

\begin{definition}
We say cylinders $C_1,\ldots, C_n$ share an angle at $P$ if $P$ is on the boundary of each and $\cap_{i=1}^n \alpha(P,C_i)\neq\emptyset$.
\end{definition}

Given a maximal cylinder $C$, let $C(\Q)$ denote the set of points in the interior of $C$ that are periodic under the action of an affine Dehn twist in $C$.  

\begin{lemma}
\label{lem:four:cylinders}
Let $C_1$ and $C_2$ are cylinders in different directions that share an angle at $P$.  Let $C_2'$ be the image of $C_2$ under a parabolic element that stabilizes the direction of $C_1$ and let $C_1'$ be the image of $C_1$ under a parabolic element that stabilizes the direction of $C_2$.  Suppose that there is a point $Q\in C_1(\Q)\cap C_1'(\Q)\cap C_2(\Q)\cap C_2'(\Q)$ and such that there is a segment from $P$ to $Q$ contained in $C_1\cap C_1'\cap C_2\cap C_2'$.  Then $C_1$ and $C_2$ are rationally related.  
\end{lemma}
\begin{proof}
Let $\alpha$ be the component of $\alpha(P,C_1)\cap\alpha(P,C_2)$ containing the initial direction of this segment.  After a normalization using $\GL(2,\R)$ we may assume $C_1$ is represented by a horizontal unit square and $C_2$ is vertical and choose local coordinates so that any $v\in\alpha$ has an initial segment in the first quadrant with $P$ at the origin.  Suppose $Q\in C_1(\Q)\cap C_1'(\Q)\cap C_2(\Q)\cap C_2'(\Q)$ lies in the first quadrant and has coordinates $(x_0,y_0)$.  
  
Since $(x_0,y_0)\in C_1(\Q)$, we have $y_0\in\Q$.  Let $x_1$ be the maximum length of a horizontal arc contained in $C_2$.  Equivalently, this is the height of the vertical cylinder $C_2$.  Since $(x_0,y_0)\in C_2(\Q)$, we have $\frac{x_0}{x_1}\in\Q$.  Since $C_1$ is a square, the inverse slope of $C_2'$ is a nonzero integer $n$ so that the line of the same slope containing $(x_0,y_0)$ is given by $x-ny=x_0-ny_0$.  The $x$ intercept is $x_0-ny_0$ and since $(x_0,y_0)\in C_2'(\Q)$, we have $$\frac{x_0-ny_0}{x_1}\in\Q,$$ which implies $\frac{y_0}{x_1}\in\Q$, and thus $x_1\in\Q$ and also $x_0\in\Q$.  Similarly,, the slope of $C_1'$ is a nonzero integer multiple $n'$ of $\frac{y_1}{x_1}$ where $y_1$ is the length of the core curve of $C_2$.  The line through $(x_0,y_0)$ of the same slope is given by $x_1y-n'xy_1=x_1y_0-n'x_0y_1$.  Since $(x_0,y_0)\in C_1'(\Q)$, the $y$-intercept $$y_0-\frac{n'x_0y_1}{x_1}\in\Q$$ so that $y_1\in\Q$.  Since $x_1\in\Q$ as well, this means $C_1$ and $C_2$ are rationally related.  
\end{proof}

\begin{definition} 
Assume that $C_1$ and $C_2$ are rationally related cylinders.  A direction is said to be \emph{rational} with respect to the pair $(C_1,C_2)$ if its slope becomes rational, after normalizing by the action of $\SL(2,\R)$ so that $C_1$ is horizontal and $C_2$ is vertical and the core curves of both cylinders are of rational length.  
\end{definition}

\begin{lemma}
\label{lem:six:cylinders}
Let $(X,\omega)$ be a non-arithmetic Veech surface.  For any singularity $P$ there are cylinders $C_i,i=1,\dots,n, n\leq 4$ each containing $P$ on its boundary such that 
\begin{itemize}
\item $\cap_{i=1}^n \alpha(P,C_i)\neq \emptyset$.
\item for any sufficiently small neighborhood $V$ of $P$,
 $$\cap_{i=1}^nC_i(\Q)\cap V=\emptyset.$$ 

\end{itemize}
\end{lemma}
\begin{proof}
 Let $C_1,C_2$ be cylinders sharing an angle at $P$.  Let $\alpha$ be any component of $\alpha(P,C_1)\cap\alpha(P,C_2)$.  By interchanging indices if necessary we may normalize so that $C_1$ is a horizontal cylinder represented by a unit square and $C_2$ is vertical.  Choose coordinates so that $P$ is at the origin and $\alpha$ is the first quadrant.  Let $\gamma_i,i=1,2$ be the generators of the parabolic subgroups stabilizing the direction of $C_i$ such that each has positive translation length.  Let $C_1'=\gamma_2(C_1)$ and $C_2'=\gamma_1^{-1}(C_2)$.  Then $$\alpha\subset\alpha(P,C_1')\cap\alpha(P,C_2').$$  We may assume $C_1$ and $C_2$ are rationally related, for otherwise we are done
by Lemma~\ref{lem:four:cylinders}.  

\textbf{Claim.} There exists a periodic direction with irrational slope.  Suppose not.  First note that the saddle connections on the boundary components of $C_1$ have rational length; for otherwise, since the height of the cylinder is $1$, there would exist a saddle connection with irrational slope crossing $C_1$.  Suppose a horizontal cylinder shares a (rational length) saddle connection with $C_1$.  Then an elementary calculation shows in fact that all of the saddle connections on its boundary have rational length and the cylinder has rational height. Since the surface is connected, it follows that every horizontal saddle connection has rational length and the heights of all of the horizontal cylinders are rational.  But then $(X,\omega)$ can be square-tiled, a contradiction.  This proves the claim.  

By the claim we can choose a periodic direction with irrational slope 
$\lambda$ with respect to $C_1$ and $C_2$.  We can assume this slope is  negative.  Then choose a cylinder $C_3$ in that direction such that $\alpha(P,C_3)$ contains $\alpha$.
  Let $\gamma_3$ be the generator of the parabolic in that direction with positive translation length, and set $C_1''=\gamma_3(C_1)$.  Set $C_3'=\gamma_1^{-1}(C_3)$.   We may assume $C_1$ and $C_3$ are rationally related, for otherwise we are again done.  

Let $V$ be any  ball of radius $\epsilon$ centered at $P$ such that 
$\epsilon$ is less than the minimum of the heights of all seven  cylinders 
$$\C_1,C_1',C_1'',C_2, C_2', C_3, C_3'.$$  
Suppose $(x_0,y_0)\in C_1(\Q)\cap C_3(\Q)\cap V$.  The equation of the line through $(x_0,y_0)$ in the direction of $C_3$ is given by $$(y-y_0)=\lambda(x-x_0).$$ 
 Because $(x_0,y_0)\in C_3(\Q)$, the $x$-intercept is rational.   Since $0\neq y_0\in \Q$, this implies that $x_0$ is irrational, so that $(x_0,y_0)\not\in C_2(\Q)$.  
\end{proof}

\begin{prop}
\label{prop:angle}
For any singularity $P$ there is a neighborhood $V$ of $P$ and a 
finite collection $\{C_i\}_{i=1}^n$ of cylinders with $P$ on their boundary such that for any $Q\in V$ there is a cylinder $C_i$ in the collection containing the segment from $P$ to $Q$ and such that $Q\notin C_i(\Q)$.
\end{prop}
\begin{proof}
We will construct the collection of cylinders inductively. By Lemma~\ref{lem:six:cylinders} there is some neighborhood $V_1$ of $P$ and $C_1^1,\ldots, C_4^1$ such that $\cap_{i=1}^4 C_1^i(\Q)\cap V_1=\emptyset$ and $\cap_{i=1}^4\alpha(P,C_i)\neq\emptyset$.  Let $\alpha_1$ be a component of the latter intersection.   Fix a small $\epsilon>0$.  Inductively suppose we have constructed an angle $\alpha_i$ at $P$, a collection of cylinders
 $C_1^i,\ldots, C_4^i$ and a neighborhood $V_i$ of $P$ such that
$$\cap_{j=1}^4 C_j^i(\Q)\cap V_i=\emptyset.$$
  Let $v_i^-, v_i^+$ be the tangent vectors that bound the angle $\alpha_i$, oriented counterclockwise. We will now construct a further collection of cylinders $C_1^{i+1},\ldots, C_4^{i+1}$.

 Let $v$ be the direction of some saddle connection inside $\alpha_i$ within $\epsilon$ of $v_i^+$.  Let $\gamma_v$ be the generator of the parabolic in the direction of $v$ with positive translation length.  For $n$ a sufficiently large negative number, $\gamma^n(v_i^-)$ is within $\epsilon$ of $v$ and the angle $\alpha_{i+1}$ between $\gamma^n(v_i^-)$ and $\gamma^n(v_i^+)$ is at least $\pi-\epsilon$.  Let $$C_j^{i+1}=\gamma^n(C_j^i), j=1,\ldots, 4,$$ and let $V_{i+1}=\gamma^n(V_i)$. Note that $$C_j^{i+1}(\Q)=\gamma^n(C_j^i(\Q)).$$  Therefore by the induction hypothesis $$\cap_{j=1}^4 C_j^{i+1}(\Q)\cap V_{i+1}=\emptyset$$ and the 
construction is complete.   

Clearly for $\epsilon$ sufficiently small, the angle conditions above imply there is $N$ depending only on the cone angle at $P$ such that every unit vector is contained in $\cup_{i=1}^N \alpha_i$.  Let $C_1,\ldots, C_n$ be given by $\cup_{i=1}^N\cup_{j=1}^4 C_j^i$.  Let $V$ be a ball of radius $r$ centered at $P$, with $r$ chosen small enough, so that $V\subset \cap_{i=1}^N V_i$, and such that $r$ is smaller than the minimum height of the cylinders $C_i$.  

Given $Q\in V$, there is a unique segment contained in $V$ joining $P$ to $Q$.  Let $i$ be any index such that $\alpha_i$ contains the initial direction of this segment.  Since $\cap_{j=1}^4 C_j^i(\Q)\cap V_i=\emptyset$, there exists a $j$ such that $Q\notin C_j^i(\Q)$.  
\end{proof}

\begin{coro}\label{cor:compact}
Given a compact set $K$ in $\SL(2,\R)/\SL(X,\omega)$ there exists a $\delta>0$ such that for any $(X,\omega)\in K$ and for any slit $\cs\subset(X,\omega)$ of length less than $\delta$ with one endpoint a zero, there exists a cylinder $C$ such that $\cs\subset C$ and such that the other endpoint of $\cs$ is not in $C(\Q)$.  
\end{coro}
\begin{proof}
For each $(X,\omega)\in K$, let $C_1,\ldots,C_n$ be the cylinders given by Proposition~\ref{prop:angle}.  Let $h$ be the minimum height of the cylinders.  Since the property of being in $C(\Q)$ is $\SL(2,\R)$-equivariant, there is a neighborhood $W$ of the identity in $\SL(2,\R)$ such that the neighborhood $A(V)$ and the collection of cylinders 
$$\{A(C_i): A\in W\}$$ satisfy the conclusion of Proposition~\ref{prop:angle} on the surface $A(X,\omega)$ and the minimum height of the cylinders is at least $h/2$.  Now cover $K$ with such open sets and take a finite subcover.  Take $\delta$ to be the minimum height of any cylinder associated with an element of the finite subcover.  
\end{proof}

\section{Proof of Theorem~\ref{thm:general}}
Let $D$ be the  $\SL(2,\R)$ orbit of 
$(X,\omega)$.  Let $K\subset D$ consist of those surfaces that have a saddle connection of length $1$. We first show $K$ is compact. 
Given a Veech surface $(X_0,\omega_0)\in K$, let $\gamma$ be a saddle connection of length one.  Since there are only a finite number of ratios of lengths of saddle connections in any cylinder decomposition, the circumferences of the cylinders in the direction of $\gamma$ have lengths bounded above.  
Since the area of the surface is one, and there are only finitely many ratios of heights of cylinders, the heights are bounded below.  It follows that the length of any saddle connection is bounded below.  This shows $K$ is a compact set.  
Now starting with the surface $(X,\omega)$ 
the idea is to apply Lemma~\ref{keylem:twist} to inductively construct sequences of slits determining small changes of area.  In order to apply that Lemma we need to know that there is a slit with small ratio so that the construction can begin.  The next Lemma says that this is possible.  

\begin{lemma}\label{small-ratio}
Let $(X,\omega)$ be a non-arithmetic Veech surface.  Let $\cs$ be a slit joining a zero $P_0$ to an aperiodic point $P_1$.  Let $(\Hat X,\Hat\omega)$ be the induced cover by $\cs$.  For any $\delta>0$ there exists a slit $\cs'$ on $X$ whose induced cover is also $(\Hat X,\Hat\omega)$ and whose ratio $\rho(\cs')<\delta$.  
\end{lemma}
\begin{proof}

Let $D$ be the $\SL(2,\R)$-orbit of $(X,\omega)$.  
Let $\MMM_D(\beta)$ be the moduli space of double covers over $D$, branched over 
  a zero and a regular point.  
Since $P_1$ is an aperiodic point, $(\Hat X,\Hat\omega)$ is not a Veech surface.  
By Theorem~\ref{thm:dense} the $\SL(2,\R)$-orbit of $(\Hat X,\Hat\omega)$ is dense 
  in the connected component of $\MMM_D(\beta)$ containing it.  

Let $K'$ be a compact neighborhood of $(X,\omega)$ in $D$.  
The lengths of closed curves and saddle connections on surfaces in $K'$ are bounded 
  below by some $\delta_0>0$.  
Each surface in $K'$ is of the form $A(X,\omega)$ for some $A$ that lies in a 
  compact neighborhood of the identity in $\SL(2,\R)$.  
Let $L$ be a set of covers in $\MMM_D(\beta)$ obtainable by applying the slit 
  construction to a surface $A(X,\omega)$ in $K'$ along a slit of length less 
  than $\delta\delta_0$ with one endpoint at $AP_0$.  
Then since the length of saddle connections of surfaces in $K'$ is at least $\delta_0$, each surface in $L$ is induced by a slit whose ratio is less than $\delta$.  
By Lemma~\ref{lem:connected}, $L$ lies in the same connected component of 
  $\MMM_D(\beta)$ that contains the given cover $(\Hat X,\Hat\omega)$.  
Moreover, since $L$ has nonempty interior, there is a surface in $L$ that 
  lies in the $\SL(2,\R)$-orbit of  $(\Hat X,\Hat\omega)$.  
Since the $\SL(2,\R)$ action respects the slit construction as well as ratios 
  of slits, $(\Hat X,\omega)$ is induced by a slit of ratio less than $\delta$.  
 \end{proof}

The following Proposition inductively constructs sequences of slits with the desired properties. 
\begin{prop}
\label{prop:children}
Let $(\hat X,\hat\omega)$ be a cover over the non-arithmetic Veech surface 
$(X,\omega)$ induced by the slit $\cs$.  For each $n$, there exist $2^n$ slits $\cs_n^j, j=1,\ldots, 2^n$ determining partitions $\Omega_n^j,\Omega_n^{j'}$ of $(\hat X,\hat \omega)$ such that for each $\cs_n^j$, called the parent, there are a pair of slits $\cs_{n+1}^i, i=1,2$ called the children, with 
\begin{itemize}
  \item the angle $\theta_n=\angle \cs_{n+1}^i \cs_n^j$ between each parent and its two children satisfies $\theta_n\leq\frac{2^{-n}}{|\cs_n^j|^2}$.
  \item Each of the $2^{n+1}$ possible angles $\theta_n$ satisfy $\theta_n\leq\delta_n/4$ where $\delta_n=\min_{i,j}\angle \cs_n^i\cs_n^j$ is the minimum angle between any two slits at level $n$. 
  \item the change of areas satisfy $a(\cs_n^j,\cs_{n+1}^i)\leq 2^{-n-1}$.
\end{itemize}
\end{prop}
\begin{proof}  
The slits are constructed inductively.  Lemma~\ref{small-ratio} allows us to find a first slit $\cs_0$ such that $$\rho(\cs_0)<\delta$$ where $\delta$ is the constant given by Corollary~\ref{cor:compact} applied to the set $K$.  
%Recall $C(l)$ is the constant in Lemma~\ref{ratio-cross-product} and $\delta$ is the constant in Corollary~\ref{cor:compact}.  Since the surface is a Veech surface, we can assume that $C(l)$ is bounded by some $C$.  Lemma~\ref{small-ratio} allows us to find a first slit $\cs_1$ such that $$\rho(\cs_1)\leq \min(\frac{1}{4C},\delta).$$  
Now suppose we have constructed $2^n$ slits $\cs_n^j$ at level $n$.  

There are two cases.  If the slit $\cs_n^j$ is on a saddle connection $\gamma$, we use the $\SL(2,\R)$-action to make its length one so that the corresponding surface is in $K$.  Since $\rho(\cs_n^j)<\delta$, Corollary~\ref{cor:compact} implies there is a cylinder $C$ containing $\cs_n^j$ so that by Lemma~\ref{keylem:twist} we can twist in a sub-cylinder of $C$ to form a pair of new slits $\cs_{n+1}^i$ with ratio $$\rho(\cs_{n+1}^i)<\min(\delta,\frac{1}{2^{n+1}C(\ell)})$$ where $C(\ell)$ is the constant in Lemma~\ref{ratio-cross-product} and $\ell$ is the maximum number of times $\gamma$ crosses a cylinder in the same direction as $C$.  
%\begin{equation}
%\label{eq:ratio}
%\rho(\cs_{n+1}^i)<\min(\frac{1}{2}\rho(\cs_n^j),\frac{1}{2^{n+1}C(\ell)}).
%\end{equation}  
By Lemma~\ref{ratio-cross-product} the change of area satisfies $$a(\cs_n^j,\cs_{n+1}^i)\leq 2^{-n-1}.$$  Since $\cs_n^j$ lies on a saddle connection $\gamma$, there exists a Dehn twist $\tau_{\cs_n^j}$ in that direction which fixes $\cs_n^j$.  Applying the twist to a power $k$, and setting $\cs_{n+1}^i=\tau_{\cs_n^j}^k(u_{n+1}^i)$, the direction of $\cs_{n+1}^i$ converges to the direction of $\cs_n^j$ as $k\to\infty$.  Thus we can choose $k$ large enough so that the first two conclusions in the Proposition hold as well.  Note that as $\tau_{\cs_n^j}^k$ is an affine map, it preserves areas. 

The second case is if $\cs_n^j$ is not on a saddle connection.  Inductively, let $\theta_n$ be the direction of $\cs_n^j$.  We apply the circular flow $R_{\pi/2 - \theta}$ so that $\cs_n^j$ is vertical and then the Teichm\"uller geodesic flow $g_t$.  As $(X,\omega)$ is a Veech surface, the linear flow in the direction $\theta$ is uniquely ergodic and therefore by Theorem~3.8 of \cite{MaTa}, the geodesic $g_t \ R_{\pi/2 - \theta}$ is recurrent.  For a sequence $t_j\to\infty$, the slit length $|g_{t_j} \ R_{\pi/2 - \theta}(\cs_n^j)|\to 0$, but all the other parameters are bounded away from $0$.  Consequently, for a sequence $t_j\to \infty$, the surface $g_{t_j}\ R_{\pi/2 - \theta}(\hat X,\hat\omega)$ contains a maximal cylinder $\CCC$ that contains $g_{t_j} \ R_{\pi/2 - \theta}(\cs_n)$ and which has bounded height, and length bounded away from $0$.  Thus there is a subcylinder $\CCC'$ which contains $g_{t_j} \ R_{\pi/2 - \theta}(\cs_n^j)$ which has area approaching $0$ as $t_j\to \infty$.  

Now, we apply a Dehn twist along $\CCC'$ an even number of times and get a pair of new slits $\cs_{n+1}^i$.  The change of area from the two consecutive partitions on $g_{t_j} \ R_{\pi/2 - \theta}(\hat X, \hat \omega)$ is bounded by $2\area(\CCC')$ by Lemma~\ref{lemma:area}.  As the construction is $\SL(2,\R)$ equivariant, we can choose $t_j$ so that the change of area $$a(\cs_n^j,\cs_{n+1}^i)\leq 2^{-n-1}.$$  Moreover it is easily checked that the angle between successive slits $\cs_n^j$ and $\cs_{n+1}^i$ is $$O\left(\frac{1}{|\cs_n^j|e^{t_j/2}}\right),$$ which implies that the first two conclusions in the Proposition hold.  
\end{proof}

To prove Theorem~\ref{thm:general} we create a binary rooted tree of slits.  At level $n$ we have $2^n$ slits, where each such slit has two children at level $n+1$.  The lengths $|\cs_n|$ of the slits are increasing and $|\cs_n|\to\infty$.  Consider any geodesic $\cs_1,\cs_2,\ldots$ of slits in the tree.  The first condition in the above Lemma says that for $i\geq n$,
 $$|\theta_i|=\angle \cs_i \cs_{i+1}\leq \frac{1}{2^i|\cs_i|^2}\leq \frac{1}{2^i|\cs_n|^2}.$$  
This shows that the sequence of angles form a Cauchy sequence and hence converge to a limiting direction.  The second condition on angles in Proposition~\ref{prop:children} implies that distinct geodesics in the tree have distinct limiting directions. Hence there are an uncountable number of limiting directions.  Since there are only countably many directions which are not minimal directions, if one can prove that each direction is not ergodic then we will have found uncountably many minimal nonergodic directions. 

To do this we need to check the conditions of Theorem~\ref{thm:criterion}.  The area condition (1) is implied by the third condition in Proposition~\ref{prop:children}.  Condition (2) is automatic since the partitions are symmetric; the components have equal area.  We check condition (3).  By a rotation we can assume the limiting direction is vertical.  Let $\alpha_n$ be the angle $\cs_n$ makes with the vertical direction.  By the first conclusion of Proposition~\ref{prop:children}, summing the infinite series, we have $$|\alpha_n|\leq \left|\sum_{i=n}^\infty \theta_i\right|\leq \frac{1}{2^{n-1}|\cs_n|^2}.$$  Now letting $h_n,v_n$ be the horizontal and vertical holomony of the slit $\cs_n$, we have $$\frac{|h_n|}{|v_n|}=\tan(\alpha_n)\leq \frac{1}{|\cs_n|^2}\leq \frac{1}{|v_n|^2}.$$  This implies $h_n\to 0$ as required.  

 Therefore, we can produce an uncountable number of limiting directions satisfying the criteria in Theorem~\ref{thm:criterion}.  This finishes the proof of Theorem~\ref{thm:general}.

We note that Theorem~\ref{thm:main} is a special case of Theorem~\ref{thm:general} applied when $P$ is a connection point.

\section{Proof of Theorem~\ref{thm:billiard}}
\begin{figure}
\begin{center}
\includegraphics{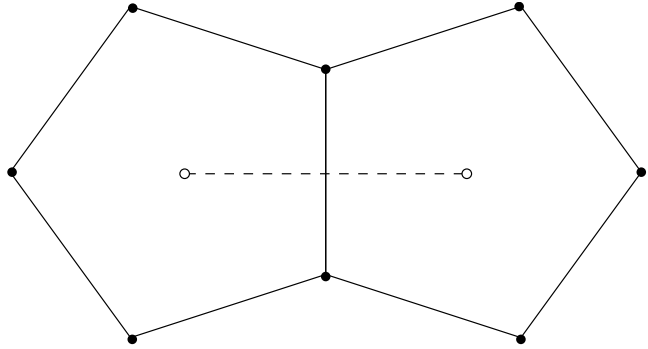}
\caption{The Veech surface $(X,\omega)$ represented as a double pentagon with opposite 
  sides identified.  The centers of the regular pentagons are joined by a horizontal slit 
  which induces the cover $(\Hat X,\Hat\omega)$ associated with the rational billiard 
  table of the $(\frac{2\pi}{5},\frac{3\pi}{10},\frac{3\pi}{10})$ isosceles triangle.}  
\label{fig:double:pentagon}
\end{center}
\end{figure}

Let $\Delta$ be the rational billiard table of the 
  $(\frac{2\pi}{5},\frac{3\pi}{10},\frac{3\pi}{10})$ triangle.  
The associated translation surface $(\Hat X,\Hat\omega)$ obtained by a standard unfolding 
  process is a branched double cover obtained by a slit construction using a slit joining 
  the centers of the double regular pentagon $(X,\omega)$.  
See Figure~\ref{fig:double:pentagon}.  

First, we prove topological dichotomy.  Let $P$ and $Q$ denote the centers of the pentagons.  The points  $P$ and $Q$ are connection points (see \cite{HuSc}). There are interchanged by the hyperelliptic involution. A saddle connection $\hat \gamma$ on $\Hat X$ projects to $\gamma$ on $X$ which is either 
\begin{enumerate}
\item \label{sc}a saddle connection, or 
\item \label{loop} a (regular) closed loop, or 
\item \label{segment} a segment joining $P$ (resp $Q$) to a zero, or 
\item \label{segmentPQ} a segment joining $P$ to $Q$ (resp. $Q$ to $P$).
\end{enumerate}

In case \ref{sc} or \ref{loop}, there is nothing to prove. The direction of $\gamma$ is periodic on $X$ as $X$ is a Veech surface; consequently the direction of $\hat \gamma$ is periodic on $\Hat X$. In case 
\ref{segment}, $\gamma$ is contained in a saddle connection as $P$ is a connection point. Thus the direction is periodic. In case \ref{segmentPQ}, the image of $\gamma$ under the hyperelliptic involution is a segment $\gamma'$ parallel to $\gamma$ joining $P$ to $Q$. Thus $\gamma \cup \gamma'$ is a closed loop on $X$ and therefore the direction of $\gamma$ is periodic. This proves that $(\Hat X,\Hat\omega)$ satisfies topological dichotomy. 

%they have coordinates in the holonomy field.  Hence any segment that joins them is either vertical or has slope in the holonomy field.  Since $(X,\omega)$ has strong holonomy type (Prop.~7.1 of \cite{HuSc}) it follows that the direction of any segment joining $P$ to $Q$ is periodic.  It follows that the direction of any saddle connection on $(\Hat X,\Hat\omega)$ induces a cylinder decomposition.  This proves that $(\Hat X,\Hat\omega)$ satisfies topological dichotomy.  

Now we show that $(\Hat X,\Hat\omega)$ does not satisfy strict ergodicity.  
By a \emph{slit} we shall mean any segment that joins $P$ and $Q$ without 
  passing through the zero and whose midpoint is a Weierstrass point.  
As noted above, the direction of any slit is periodic.  
Hence, it either lies on a saddle connection or a closed loop.  
We define the \emph{ratio} of a slit to be the ratio of its length to that 
  of the loop or saddle connection that it lies on.  

To get started, we again need to find a slit with small ratio.  
Again let $D$ be the $\SL(2,\R)$-orbit of the double pentagon $(X,\omega)$.  
Let $\MMM_D(\beta)$ be the moduli space of branched double covers over 
  surfaces in $D$ branched over $2$ regular points.  
It lies in the space $\HHH(2,1,1)$.  Let $\MMM_0\subset\MMM_D(\beta)$ denote 
  the closed $\SL(2,\R)$-invariant subset consisting of double covers branched 
  over two points that are interchanged by the hyperelliptic involution.

It was shown in \cite{EsMaMo} that the $\SL(2,\R)$-orbit of $(\Hat X,\Hat\omega)$ 
  is dense in the connected component of $\MMM_0$ that contains it.  
%Technical detail: it was shown in \cite{EsMaMo} that the $\SL(2,\R)$-orbit of 
% $(X,P,Q)$ is dense in the space $\XXX_0$ of surfaces in $D$ with a pair of 
% marked points interchanged by hyperelliptic involution.  Since the space 
% $\MMM_0$ is a unbranched cover of $\XXX_0$, dense orbits lift to dense orbits.  
Exactly the same argument as in Lemma~\ref{lem:connected} shows that every 
  branched cover in $\MMM_D(\beta)$ is given by a slit construction and if 
  two surfaces are given by a slit construction in which the slits pass through 
  Weierstrass points that correspond under the $\SL(2,\R)$ action, then they 
  lie in the same connected component of $\MMM_0$.  
Since $(\Hat X,\Hat\omega)$ can be represented as a slit construction through 
  any of the five regular Weierstrass points, $\MMM_0$ is connected.  
It follows that the $\SL(2,\R)$-orbit of $(\Hat X,\Hat\omega)$ is dense in 
  $\MMM_0$ and an argument similar to that in Lemma~\ref{small-ratio} implies 
  that we can find a slit with small ratio.  

For the induction step, given a slit $\cs$ of small ratio, we need to 
  show that there is a cylinder containing the slit so that we can apply 
  Lemma~\ref{keylem:twist} to get new slits with arbitrarily small ratio. 
To accomplish this, we will exploit properties of $(X,\omega)$ coming from 
  the fact that it belongs to the stratum $\HHH(2)$ (see $\S$  \ref{periodic-directions-inH(2)} for the definitions of the combinatorial properties of periodic directions in $\HHH(2)$).  
There are five cases.  We refer to Figure~\ref{fig:standard-picture}

The first case is if $\cs$ lies on a saddle connection $\gamma$ and 
  $\rho(\cs)\notin\Q$.  
Let $(C_1,C_2)$ be the cylinder decomposition in the direction of $\gamma$ 
  where $C_1$ is the non self-gluing cylinder and $C_2$ the self-gluing one.  
Then $\gamma$ is a saddle connection on the boundary of $C_2$ and the midpoint of $\gamma$ is a Weierstrass point. 
The closure of $C_2$ is a slit torus $T$ containing $\gamma$ in its interior.  
Lemma~\ref{keylem:twist} can now be applied to any cylinder $C\subset T$ 
  that is crossed by $\gamma$ exactly once.  

The second case is if $\cs$ lies on a saddle connection $\gamma$ and 
  $\rho(\cs)\in\Q$.  
Let $(C_1,C_2)$ be the decomposition in the direction of $\gamma$ as in the in the first case. 
Let $\gamma''$ be a saddle connection that crosses $C_2$ exactly once, joining an endpoint of $\gamma$ to itself. 
The direction of $\gamma''$ determines a cylinder decomposition and $\gamma$ crosses one of the cylinders $C''$ exactly once.  
We may Dehn twist $\gamma$ about the core curve of $C''$ exactly $3$ times to produce a  saddle connection $\gamma'$ contained in $C_2$ that crosses 
  $\gamma$ twice in its interior while missing the midpoint of $\gamma$.  
Let $(C_1',C_2')$ be the decomposition in the direction of $\gamma'$.  
Note that $\gamma$ crosses both cylinders and $\sigma$ is contained in the interior 
  of one of them.  
Let us normalize so that $\gamma$ is horizontal and $\gamma'$ is vertical.  
Let $h_i'$ be the horizontal distance across $C_i'$, i.e. its height.  
Then $|\gamma|=ah_1'+bh_2'$ for some positive integers $a$ and $b$.  
Since $(X,\omega)$ is non-arithmetic, the ratio $\frac{h_1'}{h_2'}$ is 
  irrational.  
Since $\rho(\cs)\in\Q$, the ratio of $\rho(\cs)$ to $h_i'$ is irrational 
  and Lemma~\ref{keylem:twist} can now be applied to the cylinder $C_i'$ 
  that contains the $\cs$.  

The third case is if $\cs$ lies on a loop $\gamma$ that is the core of 
  a non self-gluing cylinder and $\rho(\cs)\notin\Q$.  
Let $(C_1,C_2)$ be the decomposition in the direction of $\gamma$, where, 
  by hypothesis, $\gamma$ is the core of $C_1$.  
Let $C_1'$ be any non self-gluing cylinder disjoint from $C_1$ lying in 
  the complementary slit torus.  
Let $C_2'$ be the other cylinder in the direction of $C_1'$ and note that 
  $\gamma$ crosses only $C_2'$, perhaps multiple times.  
Lemma~\ref{keylem:twist} can be applied if the slit lies in the interior 
  of $C_2'$.  
However, it may happen that the slit crosses the boundary of $C_2'$, 
  in which case there is a unique point of intersection, namely the 
  Weierstrass point at the midpoint of the slit.  
In this situation, we apply Lemma~\ref{keylem:twist} to half of the slit, 
  then use the hyperelliptic involution to extend and get a new slit.  
It remains true that the new slit can be chosen to have arbitrarily small 
  ratio and that the area exchange can be made as small as desired.  

\begin{figure}
\begin{center}
\includegraphics{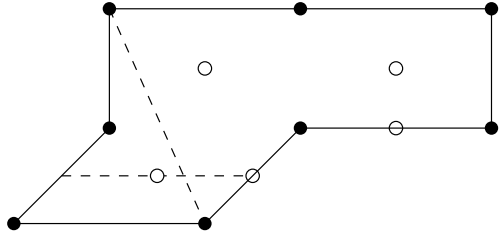}
\caption{The surface $(X,\omega)$ is normalized to have a horizontal cylinder 
decomposition $(C_1,C_2)$ where the self-gluing cylinder $C_2$ is represented 
by a rectangle such that the gluing map identifying horizontal segments on 
the boundary of $C_2$ is given by a vertical translation.  The core curve of 
$C_1$ and a saddle connection that crosses $C_1$ and $C_2$ are shown to cross 
at a point that is not a regular Weierstrass point (hollow dots).}  
\label{fig:barrel1}
\end{center}
\end{figure}

The fourth case is if $\cs$ lies on a loop $\gamma$ that is the core of 
  a non self-gluing cylinder and $\rho(\cs)\in\Q$.  
Let $(C_1,C_2)$ be the decomposition in the direction of $\gamma$, where, 
  by hypothesis, $\gamma$ is the core of $C_1$.  
We wish to find a saddle connection $\gamma'$ that is disjoint from the 
  slit $\cs$ and crosses $\gamma$ exactly once.  
Indeed, let us use the $\SL(2,\R)$ action to normalize so that $\gamma$ 
  is horizontal and $C_2$ is a represented by a rectangle such that the 
  gluing map identifying the segments on its boundary is given by a 
  vertical translation.  
Then $C_1$ may be represented by a parallelogram whose top boundary is 
  glued to the free part of the bottom boundary of $C_2$ and whose bottom boundary is glued to the free part of the top boundary of $C_2$. 
Furthermore, we may assume the parallelogram has a pair of edges with 
  positive slope such that the zero at the lower right corner lies to 
  the right of the left edge of $C_2$.  
See Figure~\ref{fig:barrel1}.  
Let $\gamma'$ be the saddle connection that joins the lower right corner 
  of the parallelogram to the zero in the interior of the top edge of 
  the rectangle representing $C_2$.  
It lies in the same homotopy class as the concatenation of the saddle 
  connection along the right edge of the parallelogram followed by the  
  saddle connection that crosses $C_2$ joining the left endpoint of self-gluing segment on the boundary of $C_2$. 
Note that $\gamma'$ does not pass through a Weierstrass point in its 
  interior and crosses $\gamma$ exactly once, as desired.  
Now, let $(C_1',C_2')$ be the decomposition in the direction of $\gamma'$.  
Note that $\gamma$ crosses both cylinders $C_i'$ since it crosses $\gamma'$.  
The rest of the argument now follows the same pattern as in the second case.  

The last case is if $\cs$ lies on a loop $\gamma$ that is the core of 
  a self-gluing cylinder.  
Let $(C_1,C_2)$ be the decomposition in the direction of $\gamma$, where, 
  by hypothesis, $\gamma$ is now the core of $C_2$.  
We normalize as in the previous case so that $C_2$ is a rectangle and the gluing  map identifying edges on its boundary is a vertical translation.   
Let $(C_1',C_2')$ be the vertical cylinder decomposition where $C_1'$ is 
  the non self gluing cylinder contained in the closure of $C_2$.  
Then $\gamma$ crosses $C_1'$ exactly once and $C_2'$ some $n'$ times.  
The length of $\gamma$ equals $h_1'+n'h_2'$ where $h_i'$ is the horizontal 
  distance across $C_i'$.  

\begin{figure}
\begin{center}
\includegraphics{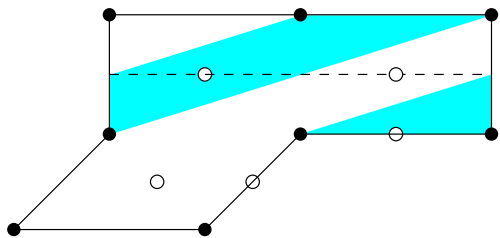}
\caption{The horizontal cylinders $(C_1,C_2)$ are represented the same way 
as in Figure~\ref{fig:barrel1}.  The vertical cylinders $(C_1',C_2')$ are 
such that the core curve of $C_2$ crosses the non self-gluing cylinder $C_1'$ 
exactly once.  The shaded region indicates the non self-gluing cylinder 
of the decomposition $(C_1'',C_2'')$ in the direction of the diagonal of $C_2$ 
of positive slope.  Note that one Weierstrass point is contained in the 
  interior of $C_1''\cap C_2'$ while the other is in $C_1'\cap C_2''$.}  
\label{fig:barrel2}
\end{center}
\end{figure}

Let $(C_1'',C_2'')$ be the decomposition in the direction of a diagonal 
  of the rectangle representing $C_2$.  
Again $\gamma$ crosses $C_1''$ once and $C_2''$ some $n''$ times.  
The length of $\gamma$ is $h_1''+n''h_2''$ where $h_i''$ is the horizontal 
  distance across $C_i''$.  
Note that $h_1'=h_1''$ because $C_1'$ and $C_1''$ are related by a Dehn 
  twist in $C_2$, not necessarily induced by an element of $\SL(X,\omega)$.  
See Figure~\ref{fig:barrel2}.  

Notice there are two Weierstrass points on $\gamma$, one of which lies 
  in $C_1''\cap C_2'$ and the other is in $C_1'\cap C_2''$.  
If the ratio of the length of $\cs$ to $h_1'$ is irrational, we can 
  apply Lemma~\ref{keylem:twist} to whichever of $C_1'$ or $C_1''$ that 
  contains $\cs$.  
Otherwise, the ratios of the length of $\cs$ to $h_2'$ and to $h_2''$ 
  are both irrational because both $\frac{h_1'}{h_2'}\notin\Q$ and 
  $\frac{h_1''}{h_2''}\notin \Q$.  
Hence, Lemma~\ref{keylem:twist} can be applied to whichever of $C_2'$ 
  or $C_2''$ that contains $\cs$.  
This completes the proof of Theorem~\ref{thm:billiard}.

\end{document}